\pdfoutput=1
\documentclass[pra,aps,twocolumn,superscriptaddress,dvipsnames,longbibliography]{revtex4}

\usepackage{amsmath,amssymb,amsthm,amsfonts,mathptmx}
\usepackage{color}
\usepackage{pbox}
\usepackage[]{algorithm2e}
\usepackage{hyperref}
\usepackage{graphicx}

\def \x {\mathbf{x}}
\def \y {\mathbf{y}}
\def \u {\mathbf{u}}
\def \v {\mathbf{v}}
\def \w {\mathbf{w}}

\begin{document}

\title{Improved stochastic trace estimation using mutually unbiased bases}

\author{J. K. Fitzsimons}
\address{Department of Engineering Science, University of Oxford}
\author{M. A. Osborne}
\address{Department of Engineering Science, University of Oxford}
\author{S. J. Roberts}
\address{Department of Engineering Science, University of Oxford}
\author{J. F. Fitzsimons}
\address{Engineering Product Development, Singapore University of Technology and Design}
\address{Centre for Quantum Technologies, National University of Singapore}


\begin{abstract}
We examine the problem of estimating the trace of a matrix $A$ when given access to an oracle which computes $\x^\dagger A \x$ for an input vector $\x$. We make use of the basis vectors from a set of mutually unbiased bases, widely studied in the field of quantum information processing, in the selection of probing vectors $\x$. This approach offers a new state of the art single shot sampling variance while requiring only $\mathcal{O}(\log(n))$ random bits to generate each vector. This significantly improves on traditional methods such as Hutchinson's and Gaussian estimators in terms of the number of random bits required and worst case sample variance.
\end{abstract}

\maketitle

\section{Introduction}

The problem of stochastic trace estimation is relevant to a range of problems from physics and applied mathematics such as electronic structure calculations \cite{bai1998computing}, seismic waveform inversion \cite{van2011seismic}, discretized parameter estimation problems with PDEs as constraints \cite{haber2012effective} and approximating the log determinant of symmetric positive semi-definite matrices \cite{boutsidis2015randomized}. Machine learning, in particular, is an example of a research domain which has many uses for stochastic trace estimation. They have been used efficiently by Generalised Cross Validation (GCV) in discretized iterative methods for fitting Laplacian smoothing splines to very large datasets \cite{hutchinson1990stochastic}, computing the number of triangles in a graph \cite{atallah2001randomized, atallah2013lower}, string pattern matching \cite{avron2010counting,tsourakakis2008fast} and the training Gaussian Processes using score functions \cite{stein2013stochastic}.

Stochastic trace estimation endeavours to choose $n$ dimensional vectors $x$ such that the expectation of $x^TAx$ is equal to the trace of the implicit symmetrical positive semi definite matrix $A \in \mathbb{R}^{n\times n}$. It can be seen that many sampling policies satisfy this condition.  As such several metrics are used in order to choose a sampling policy such as the one sample variance, the number of samples to achieve a $(\epsilon, \delta)$-approximation and the number of random bits required to create $x$ \cite{avron2011randomized}. This last metric is motivated in part by the relatively long timescales for hardware number generation, and concerns about parallelising pseudo-random number generators.

In this work we propose a new stochastic trace estimator based on mutually unbiased bases (MUBs) \cite{schwinger1960unitary}, and quantify the single shot sampling variance of the proposed MUBs sampling method and its corresponding required number of random bits. We will refer to methods which sample from a fixed set of basis functions as being fixed basis sampling methods. For example, we can randomly sample the diagonal values of the matrix $A$ by sampling $z$ from the set of columns which form the identity matrix. This is referred to as the unit vector estimator in the literature \cite{avron2011randomized}. Other similar methods sample from the columns Discrete Fourier Transform (DFT), the Discrete Hartley Transform (DHT), the Discrete Cosine Transform (DCT) or a Hadamard matrix. We prove that sampling from the set of mutually unbiased bases significantly reduces this single shot sample variance, in particular in the worst case bound.

The paper is laid out as follows: Section \ref{MUBS} gives a brief introduction to mutually unbiased basis, Section \ref{TA} describes our novel approach of using mutually unbiased bases for trace estimation and Section \ref{AN} gives a rigorous analysis of of the new estimator. Section \ref{exp} compares the proposed MUBs estimator to established approaches both in terms of the analytic expectation of sample variance and as applied to synthetic and real data. The task of counting the number of triangles in a graph is considered as an example application.

\section{Mutually Unbiased Basis}\label{MUBS}

Linear algebra has found application in a diverse range of fields, with each field drawing from a common set of tools. However, occasionally, techniques developed in one field do not become well known outside of that community, despite the potential for wider use. In this work, we will make extensive use of mutually unbiased bases, sets of bases that arise from physical considerations in the context of quantum mechanics \cite{schwinger1960unitary} and which have been extensively exploited within the quantum information community \cite{durt2010mutually}. In quantum mechanics, physical states are represented as vectors in a complex vector space, and the simplest form of measurement projects the state onto one of the vectors from some fixed orthonormal basis for the space, with the probability for a particular outcome given by the square of the length of the projection onto the corresponding basis vector \footnote{For a more comprehensive introduction to the mathematics of quantum mechanics in finite-dimensional systems, we refer the reader to \cite{nielsen2010quantum}}. In such a setting, it is natural to ask about the existence of pairs or sets of measurements where the outcome of one measurement reveals nothing about the outcome of another measurement, and effectively erases any information about the outcome had the alternate measurement instead been performed. As each measurement corresponds to a particular basis, such a requirement implies that the absolute value of the overlap between pairs of vectors drawn from bases corresponding to different measurements be constant. This leads directly to the concept of mutually unbiased bases (MUBs).

A set of orthonormal bases $\{B_1,\ldots,B_n\}$ are said to be mutually unbiased if for all choices of $i$ and $j$, such that $i\neq j$, and for every $\u \in B_i$ and every $\v \in B_j$, $|\u^\dagger \v| = \frac{1}{\sqrt{n}}$, where $n$ is the dimension of the space. While for real vector spaces the number of mutually unbiased bases has a complicated relationship with the dimensionality \cite{boykin2005real}, for complex vector spaces the number of mutually unbiased bases is known to be exactly $n+1$ when $n$ is either a prime or an integer power of a prime \cite{klappenecker2004constructions}. Furthermore, a number of constructions are known for constructing such bases \cite{klappenecker2004constructions}. When $n$ is neither prime nor a power of a prime, the number of mutually unbiased bases remains open, even for the case of $n=6$ \cite{butterley2007numerical}, but is known to be at least $p_1^{d_1} + 1$, where $n = \prod_{i} p_i^{d_i}$ and $p_i$ are prime numbers such that $p_i < p_{i+1}$ for all $i$.  

\section{Trace Estimators}\label{TA}

In order to estimate the trace of a $n\times n$ positive semi-definite matrix $A$ from a single call to an oracle for $\x^\dagger A \x$, we consider four strategies:
\begin{itemize}
\item \textbf{Fixed basis estimator:} For a fixed orthonormal basis $B$, choose $\x$ uniformly at random from the elements of $B$. The trace is then estimated to be $n \x^\dagger A \x$.
\item \textbf{Mutually unbiased bases (MUBs) estimator:} For a fixed choice of a set of $b$ mutually unbiased bases $\mathbb{B} = \{B_1,...,B_b\}$, choose $B$ uniformly at random from $\mathbb{B}$ and then choose $\x$ uniformly at random from the elements of $B$. Here $b$ is taken to be the maximum number of mutually unbiased bases for a complex vector space of dimension $n$. As in the fixed basis strategy, the trace is then estimated to be $n \x^\dagger A \x$.
\item \textbf{Hutchinson's estimator:} Randomly choose the elements of $\x$ independently and identically distributed from a Rademacher distribution $\left(Pr(x_i = \pm 1) = \frac{1}{2}\right)$. The trace is then estimated to be $\x^\dagger A \x$.
\item \textbf{Gaussian estimator:} Randomly choose the elements of $\x$ independently and identically distributed from a zero mean unit variance Gaussian distribution. The trace is then estimated to be $\x^\dagger A \x$.
\end{itemize}
The first strategy is a generic formulation of approaches which sample vectors from a fixed orthogonal basis, the most efficient sampling method in terms of the number of random bits required in the literature \cite{avron2011randomized}, while the second strategy is novel and represents our main contribution. Both strategies have similar randomness requirements: In the first strategy at least $\lceil\log_2 (n)\rceil$ random bits are necessary to ensure the possibility of choosing every element of $B$. In the second strategy, an identical number of random bits is necessary to choose $\x$ for a fixed $B$, and $\lceil\log_2 (b)\rceil$ random bits are necessary to choose $B$. Note that an upper bound on the number of mutually unbiased bases is one greater than dimensionality of the space, and this bound is saturated for spaces where the dimensionality is prime or an integer power of a prime, i.e. $b\leq n+1$. Thus the number of random bits necessary to implement these strategies differs by a factor of approximately two. The third and forth strategies significantly outperform the fixed basis estimator in terms of single-shot variance, at the cost of a dramatic increase in the amount of randomness required, and have been extensively studied in the literature \cite{avron2011randomized, hutchinson1989stochastic, silver1997calculation}. For conciseness we will not repeat the analysis of these methods in this paper but will compare the fixed basis estimator and MUBs estimator to them in Table \ref{results}. 

\begin{table*}[t]\label{results}
\centering
\begin{tabular}{|l|c|c|c|}
\hline
Estimator~~~&~~~ $V$ ~~~&~~~ $V^\text{worst}$ ~~~&~~~ $R$ \\ \hline
&&&\\
Fixed basis ~~~&~~~  $n \sum_{i=1}^n M_{ii}^2 - \text{Tr}(A)^2$   ~~~&~~~    $(n-1)\text{Tr}(A)^2$              ~~~&~~~        $\log_2(n)$             \\
&&&\\
MUBs         ~~~&~~~ $\frac{n}{n+1} \text{Tr}(A^2) - \frac{1}{n+1} \text{Tr}(A)^2$    ~~~&~~~        $\frac{n}{n+1}\text{Tr}(A^2)$          ~~~&~~~          $\log_2(n)+\log_2(n+1)$           \\
&&&\\
Hutchinson \cite{hutchinson1989stochastic} ~~~ &~~~ $2\left(\text{Tr}(A^2) - \sum_{i=1}^n A_{ii}^2\right)$   ~~~ & ~~~      $\frac{2(n-1)}{n} \text{Tr}(A^2)$           ~~~&~~~        $n$             \\
&&&\\
Gaussian \cite{silver1997calculation} ~~~&~~~ $2\text{Tr}(A^2)$ ~~~&~~~ $2\text{Tr}(A^2)$ ~~~&~~~  $\infty$ for exact; $\mathcal{O}(n)$ for fixed precision~~~\\
\hline
\end{tabular}
\centering
\caption{Comparison of single shot variance $V$, worst case single shot variance $V^\text{worst}$ and number of random bits $R$ required for commonly used trace estimators and the MUBs estimator. In the case of the MUBs estimator, the quantities provided for the variances are upper bounds rather than the exact variance.}\label{tab:comp}
\end{table*}

\subsection{Analysis of fixed basis estimator}
We first analyse the worst case variance of the fixed base estimator. In this analysis and the analysis for the MUBs estimator which follows, we make no assumption on $A$ and consider the worst case variance.

We begin from the definition of the variance of the estimator for a single query. Let $X$ be a random variable such that $X = \x^\dagger A \x$, where $\x$ is chosen according to the fixed basis strategy. Then
\begin{equation}
\text{Var(X)} = E(X^2) - E(X)^2,\label{eq:var1}
\end{equation}
where $E(\cdot)$ denotes the expectation value of the argument. We compute this term by term. First
\begin{eqnarray}
E(X) = \frac{1}{n}\sum_{\x \in B} \x^\dagger A \x = \frac{\text{Tr}(A)}{n}.\label{eq:Ex1}
\end{eqnarray}
where $n=\dim (A)$, and hence the second term in Eq.~\ref{eq:var1} is equal to $\frac{\text{Tr}(A)^2}{n^2}$. Turning to the first term,
\begin{eqnarray}
E(X^2) &=& \frac{1}{n}\sum_{\x \in B} \left(\x^\dagger A \x\right)^2 \label{eq:branchtoworst}\\
&=& \frac{1}{n}\sum_{i=1}^n M_{ii}^2,
\end{eqnarray}
where $M = U A U^\dagger$ for some fixed unitary matrix $U$, such that $U^\dagger \x$ is a vector in the standard basis for all $x\in B$, and $M_{ii}$ is the $i$th entry on the main diagonal of $M$. The variance for the fixed basis estimator is then given by $V_\text{fixed} = n \sum_{i=1}^n M_{ii}^2 - \text{Tr}(A)^2$. The worst case occurs when the value of $\sum_{i=1}^n M_{ii}^2$ is maximized for fixed trace of $A$ (and hence $M$), which occurs when a single diagonal entry is non-zero, and so the worst case single shot variance for the fixed basis estimator is $V_{fixed}^\text{worst} = (n-1)\text{Tr}(A)^2$.

\subsection{Analysis of MUBs estimator}\label{AN}
We now turn to analysis of the MUBs estimator. We assume that $n$ is either prime or a prime raised to some integer power, since when this is not the case the matrix can always be padded out with zeros to such a dimension with little overhead. In this case, it has been established that $b=n+1$ \cite{klappenecker2004constructions}. The variance is remains as given in Eq.~\ref{eq:var1}, except that $X$ is defined in terms of vectors $\x$ chosen according to the MUBs strategy. Again, we analyse the individual terms making up the variance. We begin with
\begin{eqnarray}
E(X) = \frac{1}{nb} \sum_{B\in \mathbb{B}} \sum_{\x \in B} \x^\dagger A \x = \frac{\text{Tr}(A)}{n}.\label{eq:Ex1}
\end{eqnarray}
and hence the second term in the variance is the same as for the fixed basis estimator. Analysing the first term is, however, more difficult. We begin with the observation that $E(X^2)$ can be expressed in terms of the trace of the Kronecker product of two matrices, as follows
\begin{eqnarray}
E(X^2) &=& \frac{1}{nb}\sum_{B\in \mathbb{B}} \sum_{\x \in B} \left(\x^\dagger A \x\right)^2 \\
&=& \frac{1}{nb}\sum_{B\in \mathbb{B}} \sum_{\x \in B} \text{Tr}\left((\x \x^\dagger A)^{\otimes 2} \right).
\end{eqnarray}
Moving the summations inside the equation we obtain
\begin{eqnarray}
E(X^2) &=& \frac{1}{nb} \text{Tr}\left(\sum_{B\in \mathbb{B}} \sum_{\x \in B}\left(\x \x^\dagger\right)^{\otimes 2}  A^{\otimes 2} \right)\\
&=& \frac{2}{nb} \text{Tr}\left(P  A^{\otimes 2} \right),\label{eq:almostvar}
\end{eqnarray}
where $P=\frac{1}{2}\sum_{B\in \mathbb{B}} \sum_{\x \in B}\left(\x \x^\dagger\right)^{\otimes 2}$. 

While this form of $P$ may appear intimidating, we now prove that $P$ is in fact a projector with each eigenvalue being either $0$ or $1$. We prove this indirectly, first by showing that $P$ has rank at most $n(n+1)/2$, and then using the relationship between the traces of $P$ and $P^2$ to conclude that the remaining $n(n+1)/2$ eigenvalues are equal to unity. Any vector of the form $\w = \u\otimes \v - \u\otimes \v$ for $\u,\v\in B_1$ trivially satisfies $P \w = \textbf{0}$. Since such vectors form a basis for a subspace of dimension $n(n-1)/2$, we conclude that $\text{rank}(P) \leq n^2 - n(n-1)/2 = n(n+1)/2$. Turning now to the issue of trace, we have
\begin{eqnarray}
\text{Tr}(P) &=& \text{Tr}\left(\frac{1}{2}\sum_{B\in \mathbb{B}} \sum_{\x \in B}\left(\x \x^\dagger\right)^{\otimes 2}\right)\\
&=& \frac{1}{2}\sum_{B\in \mathbb{B}} \sum_{\x \in B}\left(\x^\dagger \x\right)^2\\
&=& \frac{nb}{2}.
\end{eqnarray}
We can similarly compute the trace of $P^2$ to obtain
\begin{eqnarray}
\text{Tr}(P^2) &=& \text{Tr}\left(\frac{1}{4}\sum_{B,B'\in \mathbb{B}} \sum_{\x \in B} \sum_{\y \in B'}\left(\x \x^\dagger\right)^{\otimes 2}\left(\y \y^\dagger\right)^{\otimes 2}\right)\\
&=& \frac{1}{4}\sum_{B,B'\in \mathbb{B}} \sum_{\x \in B} \sum_{\y \in B'} \left|\x^\dagger \y\right|^4\\
&=& \frac{nb}{4} + \frac{n^2b(b-1)}{4 n^2}\\
&=& \frac{b(n+b-1)}{4}.
\end{eqnarray}
Notice that this implies that $\text{Tr}(P) = \text{Tr}(P^2)$ for dimensions which are prime or integer powers of a prime, since in such cases $b = n+1$. This implies that the eigenvalues on the non-zero subspace minimize the sum of their squares for a fixed sum, and since $P$ is positive semi-definite, we can conclude that each non-zero eigenvalue must be equal to unity.

Returning to the calculation of variance, we then have
\begin{eqnarray}
E(X^2) &\leq& \frac{2}{nb} \text{Tr}\left(A^{\otimes 2} \right)\\
&=& \frac{2}{nb} \text{Tr}\left(A\right)^2,
\end{eqnarray}
and hence
\begin{equation}
\text{Var}(X) \leq \left(\frac{2}{nb} - \frac{1}{n^2}\right) \text{Tr}\left(M\right)^2 \leq\frac{\text{Tr}\left(A\right)^2}{n^2}.\label{eq:simplevar}
\end{equation}
This implies that the variance on the estimate of $\text{Tr}(A)$ is bounded from above by $\text{Tr}(A)^2$. It is, in fact, possible to compute the variance exactly from Eq.~\ref{eq:almostvar} by observing that $M$ is the projector onto the symmetric subspace when $n$ is an integer power of a prime. That is to say, for any vector $\u$ and any vector $\v$ orthogonal to $\u$, the vectors $\u\otimes\v + \v\otimes\u$, $\u\otimes\u$ and $\v \otimes \v$ are in the $+1$ eigenspace of $M$, whereas the vector $\u\otimes\v - \v\otimes\u$ is in the null space of $M$. We can then compute the exact variance of the MUBs estimator using the spectral decomposition $A = \sum_i \lambda_i \u_i \u_i^\dagger$ as
\begin{eqnarray}
\text{V}_\text{MUBs} &=&  \frac{2n}{n+1} \text{Tr}\left(P  A^{\otimes 2} \right) - \text{Tr}(A)^2\\
&=&\frac{2n}{n+1} \sum_{i=1}^n \sum_{j=1}^n \lambda_i \lambda_j \text{Tr}\left(P (\u_i\otimes\u_j)(\u_i\otimes\u_j)^\dagger \right) - \text{Tr}(A)^2 \nonumber\\
&=&\frac{2n}{n+1} \sum_{i=1}^n \left( \lambda_i^2 + \frac{1}{2}\sum_{j\neq i}\lambda_i \lambda_j \right) - \text{Tr}(A)^2\\
&=& \frac{n}{n+1} \text{Tr}(A^2) - \frac{1}{n+1} \text{Tr}(A)^2.
\end{eqnarray}
Since for all positive semi-definite matrices $A$ the value of $\text{Tr}(A)^2$ is bounded from below by $\text{Tr}(A^2)$, the single shot variance on the MUBs estimator is bounded by $\text{V}_\text{MUBs}^\text{worst} = \frac{n-1}{n+1} \text{Tr}(A^2)$ in the worst case, a significant improvement on the bound stemming from Eq.~\ref{eq:simplevar}. Even if when not restricted to positive semi-definite $A$, the worst case variance is bounded by $\text{V}_\text{MUBs}^\text{worst} = \frac{n}{n+1} \text{Tr}(A^2)$, since $\text{Tr}(A)^2$ is non-negative for any $A$ defined over $\mathbb{R}^n$. The worst case single shot variance of the MUBs estimator is then at least a factor of $n-1$ better than that of any fixed basis estimator. Furthermore, the variance for the widely used Hutchinson estimator \cite{hutchinson1989stochastic, avron2011randomized}, is given by $V_\text{H} = 2\left(\text{Tr}(A^2) - \sum_{i=1}^n A_{ii}^2\right)$. In the worst case, $\sum_{i=1}^n A_{ii}^2 = \frac{1}{n}\text{Tr}(A^2)$, and hence the worst case single shot variance for Hutchinson estimator is $V_\text{H}^\text{worst} = \frac{2(n-1)}{n} \text{Tr}(A^2)$. Thus, the MUBs estimator has better worst case performance than the Hutchinson estimator by a factor $\frac{2(n+1)(n-1)}{n^2}$ which approaches $2$ for large $n$. This improvement is perhaps unsurprising, since for a symmetric matrices $\x^\dagger A \x = \x_R^T A \x_R + \x_I^T A \x_I$ where $\x_R$ and $\x_I$ are the real and imaginary parts of $\x$. Hence evaluating $\x^\dagger A \x$ for a single complex vector is equivalent to taking the sum of it over for two different real vectors, leading to a factor of two improvement in the variance of the average.

Table \ref{tab:comp} compares the single shot variance, worst case single shot variance and randomness requirements of the trace estimators. As can be seen from the comparison the MUBs estimator has strictly smaller variance than either the Hutchinson or Gaussian methods, while requiring significantly less randomness to implement. Given the drastic reduction in randomness requirements, and the improved worst case performance, the MUBs estimator provides an attractive alternative to previous methods for estimating the trace of implicit matrices.

\section{Numerical Results}\label{exp}

\begin{figure}
\includegraphics[width=0.9\columnwidth]{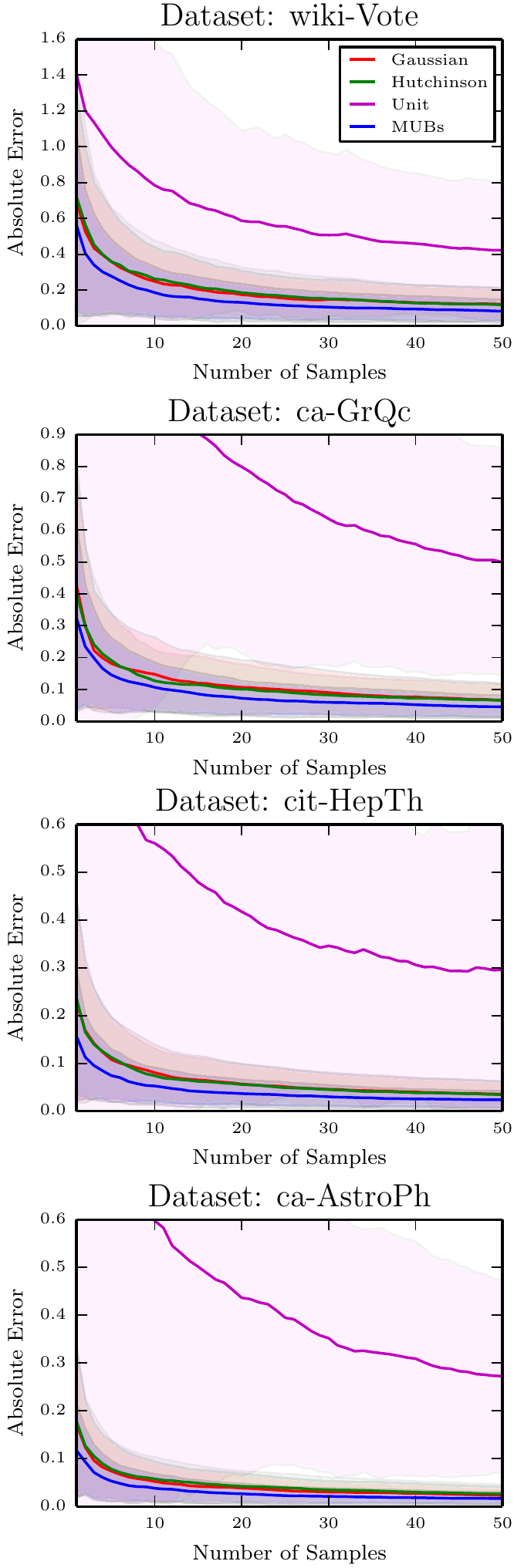}
\caption{A comparison of the performance of the stochastic trace estimation methods on the four datasets. The fixed basis method was not included as it was not competitive. The experiments were performed 500 times each. The solid line indicated the empirical mean absolute relative error and the the surrounding transparent region indicates one empirical standard deviation of the 500 trials.}\label{triangle_res}
\end{figure}

Having quantified the worst case performance of the MUBs estimator, we now explore its performance in practice, in a set of numerical experiments. As an example application we will consider counting the number of triangles in a graph. This is an important problem in a number of application domains such as identifying the number of mutual acquaintences in a social network. An efficient method to do this is the Trace Triangle algorithm \cite{avron2010counting}. The algorithm is based on a relationship between the adjacency matrix, $A$, and the number of triangles for an undirected graph, $\Delta_g$,

\begin{equation}
\Delta_g = \frac{Tr(A^3)}{6}.
\end{equation}  

The trace of the adjacency matrix cubed can be sampled in $\mathcal{O}(n^2)$ per sample as opposed to being explicitly computed in $\mathcal{O}(n^3)$. We compared Gaussian, Hutchinson's, Unit and MUBs estimators performance at predicting the number of triangles for the graphs presented in Table \ref{datasets} and the results of the experiment are presented in Figure \ref{triangle_res}. The code for these experiments, with an efficient Python implementation for generating the MUBs sample vectors in $\mathcal{O}(n)$, has been made publically available \footnote{\url{www.github.com/OxfordMLRG/traceEst}}. The MUBs estimator outperforms the classical method in all of the experiments, as would be expected from the theoretical analysis in terms of variance. In addition, the exponential reduction in randomness means that implementations making use of hardware random number generation will generally see a significant decrease in processing times.

\begin{table}[t!]\label{datasets}
\begin{tabular}{|l|l|l|l|}
\hline
Dataset &Vertices & Edges & Triangles\\
\hline
Arxiv-HEP-th & 27,240 & 341,923 &1,478,735 \\
CA-AstroPh & 18,772 & 198,050 & 1,351,441 \\
CA-GrQc & 5,242 & 14,484 & 48,260 \\
wiki-vote & 7,115 & 100,689 & 608,389 \\
\hline
\end{tabular}
\centering
\caption{Datasets used for the comparison of stochastic trace estimation methods in the counting of triangles in graphs. All datasets can be found at \url{snap.stanford.edu/data}\label{datasets}}
\end{table}

\section*{Acknowledgements}

JFF acknowledges support from the Air Force Office of Scientific Research under AOARD grant FA2386-15-1-4082. This material is based on research funded in part by the Singapore National Research Foundation under NRF Award NRF-NRFF2013-01.

\bibliographystyle{apsrev}
\bibliography{TraceEst_MUBS}

\begin{thebibliography}{19}
\expandafter\ifx\csname natexlab\endcsname\relax\def\natexlab#1{#1}\fi
\expandafter\ifx\csname bibnamefont\endcsname\relax
  \def\bibnamefont#1{#1}\fi
\expandafter\ifx\csname bibfnamefont\endcsname\relax
  \def\bibfnamefont#1{#1}\fi
\expandafter\ifx\csname citenamefont\endcsname\relax
  \def\citenamefont#1{#1}\fi
\expandafter\ifx\csname url\endcsname\relax
  \def\url#1{\texttt{#1}}\fi
\expandafter\ifx\csname urlprefix\endcsname\relax\def\urlprefix{URL }\fi
\providecommand{\bibinfo}[2]{#2}
\providecommand{\eprint}[2][]{\url{#2}}

\bibitem[{\citenamefont{Bai et~al.}(1998)\citenamefont{Bai, Fahey, Golub,
  Menon, and Richter}}]{bai1998computing}
\bibinfo{author}{\bibfnamefont{Z.}~\bibnamefont{Bai}},
  \bibinfo{author}{\bibfnamefont{M.}~\bibnamefont{Fahey}},
  \bibinfo{author}{\bibfnamefont{G.~H.} \bibnamefont{Golub}},
  \bibinfo{author}{\bibfnamefont{M.}~\bibnamefont{Menon}}, \bibnamefont{and}
  \bibinfo{author}{\bibfnamefont{E.}~\bibnamefont{Richter}},
  \bibinfo{type}{Tech. Rep.}, \bibinfo{institution}{Citeseer}
  (\bibinfo{year}{1998}).

\bibitem[{\citenamefont{van Leeuwen et~al.}(2011)\citenamefont{van Leeuwen,
  Aravkin, and Herrmann}}]{van2011seismic}
\bibinfo{author}{\bibfnamefont{T.}~\bibnamefont{van Leeuwen}},
  \bibinfo{author}{\bibfnamefont{A.~Y.} \bibnamefont{Aravkin}},
  \bibnamefont{and} \bibinfo{author}{\bibfnamefont{F.~J.}
  \bibnamefont{Herrmann}}, \bibinfo{journal}{International Journal of
  Geophysics} \textbf{\bibinfo{volume}{2011}} (\bibinfo{year}{2011}).

\bibitem[{\citenamefont{Haber et~al.}(2012)\citenamefont{Haber, Chung, and
  Herrmann}}]{haber2012effective}
\bibinfo{author}{\bibfnamefont{E.}~\bibnamefont{Haber}},
  \bibinfo{author}{\bibfnamefont{M.}~\bibnamefont{Chung}}, \bibnamefont{and}
  \bibinfo{author}{\bibfnamefont{F.}~\bibnamefont{Herrmann}},
  \bibinfo{journal}{SIAM Journal on Optimization}
  \textbf{\bibinfo{volume}{22}}, \bibinfo{pages}{739} (\bibinfo{year}{2012}).

\bibitem[{\citenamefont{Boutsidis et~al.}(2015)\citenamefont{Boutsidis,
  Drineas, Kambadur, and Zouzias}}]{boutsidis2015randomized}
\bibinfo{author}{\bibfnamefont{C.}~\bibnamefont{Boutsidis}},
  \bibinfo{author}{\bibfnamefont{P.}~\bibnamefont{Drineas}},
  \bibinfo{author}{\bibfnamefont{P.}~\bibnamefont{Kambadur}}, \bibnamefont{and}
  \bibinfo{author}{\bibfnamefont{A.}~\bibnamefont{Zouzias}},
  \bibinfo{journal}{arXiv preprint arXiv:1503.00374}  (\bibinfo{year}{2015}).

\bibitem[{\citenamefont{Hutchinson}(1990)}]{hutchinson1990stochastic}
\bibinfo{author}{\bibfnamefont{M.~F.} \bibnamefont{Hutchinson}},
  \bibinfo{journal}{Communications in Statistics-Simulation and Computation}
  \textbf{\bibinfo{volume}{19}}, \bibinfo{pages}{433} (\bibinfo{year}{1990}).

\bibitem[{\citenamefont{Atallah et~al.}(2001)\citenamefont{Atallah, Chyzak, and
  Dumas}}]{atallah2001randomized}
\bibinfo{author}{\bibfnamefont{M.~J.} \bibnamefont{Atallah}},
  \bibinfo{author}{\bibfnamefont{F.}~\bibnamefont{Chyzak}}, \bibnamefont{and}
  \bibinfo{author}{\bibfnamefont{P.}~\bibnamefont{Dumas}},
  \bibinfo{journal}{Algorithmica} \textbf{\bibinfo{volume}{29}},
  \bibinfo{pages}{468} (\bibinfo{year}{2001}).

\bibitem[{\citenamefont{Atallah et~al.}(2013)\citenamefont{Atallah, Grigorescu,
  and Wu}}]{atallah2013lower}
\bibinfo{author}{\bibfnamefont{M.~J.} \bibnamefont{Atallah}},
  \bibinfo{author}{\bibfnamefont{E.}~\bibnamefont{Grigorescu}},
  \bibnamefont{and} \bibinfo{author}{\bibfnamefont{Y.}~\bibnamefont{Wu}},
  \bibinfo{journal}{Information Processing Letters}
  \textbf{\bibinfo{volume}{113}}, \bibinfo{pages}{690} (\bibinfo{year}{2013}).

\bibitem[{\citenamefont{Avron}(2010)}]{avron2010counting}
\bibinfo{author}{\bibfnamefont{H.}~\bibnamefont{Avron}}, in
  \emph{\bibinfo{booktitle}{Workshop on Large-scale Data Mining: Theory and
  Applications}} (\bibinfo{year}{2010}), vol.~\bibinfo{volume}{10}, pp.
  \bibinfo{pages}{10--9}.

\bibitem[{\citenamefont{Tsourakakis}(2008)}]{tsourakakis2008fast}
\bibinfo{author}{\bibfnamefont{C.~E.} \bibnamefont{Tsourakakis}}, in
  \emph{\bibinfo{booktitle}{Data Mining, 2008. ICDM'08. Eighth IEEE
  International Conference on}} (\bibinfo{organization}{IEEE},
  \bibinfo{year}{2008}), pp. \bibinfo{pages}{608--617}.

\bibitem[{\citenamefont{Stein et~al.}(2013)\citenamefont{Stein, Chen, Anitescu
  et~al.}}]{stein2013stochastic}
\bibinfo{author}{\bibfnamefont{M.~L.} \bibnamefont{Stein}},
  \bibinfo{author}{\bibfnamefont{J.}~\bibnamefont{Chen}},
  \bibinfo{author}{\bibfnamefont{M.}~\bibnamefont{Anitescu}},
  \bibnamefont{et~al.}, \bibinfo{journal}{The Annals of Applied Statistics}
  \textbf{\bibinfo{volume}{7}}, \bibinfo{pages}{1162} (\bibinfo{year}{2013}).

\bibitem[{\citenamefont{Avron and Toledo}(2011)}]{avron2011randomized}
\bibinfo{author}{\bibfnamefont{H.}~\bibnamefont{Avron}} \bibnamefont{and}
  \bibinfo{author}{\bibfnamefont{S.}~\bibnamefont{Toledo}},
  \bibinfo{journal}{Journal of the ACM (JACM)} \textbf{\bibinfo{volume}{58}},
  \bibinfo{pages}{8} (\bibinfo{year}{2011}).

\bibitem[{\citenamefont{Schwinger}(1960)}]{schwinger1960unitary}
\bibinfo{author}{\bibfnamefont{J.}~\bibnamefont{Schwinger}},
  \bibinfo{journal}{Proceedings of the National Academy of Sciences}
  \textbf{\bibinfo{volume}{46}}, \bibinfo{pages}{570} (\bibinfo{year}{1960}).

\bibitem[{\citenamefont{Durt et~al.}(2010)\citenamefont{Durt, Englert,
  Bengtsson, and {\.Z}yczkowski}}]{durt2010mutually}
\bibinfo{author}{\bibfnamefont{T.}~\bibnamefont{Durt}},
  \bibinfo{author}{\bibfnamefont{B.-G.} \bibnamefont{Englert}},
  \bibinfo{author}{\bibfnamefont{I.}~\bibnamefont{Bengtsson}},
  \bibnamefont{and}
  \bibinfo{author}{\bibfnamefont{K.}~\bibnamefont{{\.Z}yczkowski}},
  \bibinfo{journal}{International journal of quantum information}
  \textbf{\bibinfo{volume}{8}}, \bibinfo{pages}{535} (\bibinfo{year}{2010}).

\bibitem[{\citenamefont{Boykin et~al.}(2005)\citenamefont{Boykin, Sitharam,
  Tarifi, and Wocjan}}]{boykin2005real}
\bibinfo{author}{\bibfnamefont{P.~O.} \bibnamefont{Boykin}},
  \bibinfo{author}{\bibfnamefont{M.}~\bibnamefont{Sitharam}},
  \bibinfo{author}{\bibfnamefont{M.}~\bibnamefont{Tarifi}}, \bibnamefont{and}
  \bibinfo{author}{\bibfnamefont{P.}~\bibnamefont{Wocjan}},
  \bibinfo{journal}{arXiv preprint quant-ph/0502024}  (\bibinfo{year}{2005}).

\bibitem[{\citenamefont{Klappenecker and
  R{\"o}tteler}(2004)}]{klappenecker2004constructions}
\bibinfo{author}{\bibfnamefont{A.}~\bibnamefont{Klappenecker}}
  \bibnamefont{and}
  \bibinfo{author}{\bibfnamefont{M.}~\bibnamefont{R{\"o}tteler}}, in
  \emph{\bibinfo{booktitle}{Finite fields and applications}}
  (\bibinfo{publisher}{Springer}, \bibinfo{year}{2004}), pp.
  \bibinfo{pages}{137--144}.

\bibitem[{\citenamefont{Butterley and Hall}(2007)}]{butterley2007numerical}
\bibinfo{author}{\bibfnamefont{P.}~\bibnamefont{Butterley}} \bibnamefont{and}
  \bibinfo{author}{\bibfnamefont{W.}~\bibnamefont{Hall}},
  \bibinfo{journal}{Physics Letters A} \textbf{\bibinfo{volume}{369}},
  \bibinfo{pages}{5} (\bibinfo{year}{2007}).

\bibitem[{\citenamefont{Hutchinson}(1989)}]{hutchinson1989stochastic}
\bibinfo{author}{\bibfnamefont{M.}~\bibnamefont{Hutchinson}},
  \bibinfo{journal}{Communications in Statistics-Simulation and Computation}
  \textbf{\bibinfo{volume}{18}}, \bibinfo{pages}{1059} (\bibinfo{year}{1989}).

\bibitem[{\citenamefont{Silver and R{\"o}der}(1997)}]{silver1997calculation}
\bibinfo{author}{\bibfnamefont{R.}~\bibnamefont{Silver}} \bibnamefont{and}
  \bibinfo{author}{\bibfnamefont{H.}~\bibnamefont{R{\"o}der}},
  \bibinfo{journal}{Physical Review E} \textbf{\bibinfo{volume}{56}},
  \bibinfo{pages}{4822} (\bibinfo{year}{1997}).

\bibitem[{\citenamefont{Nielsen and Chuang}(2010)}]{nielsen2010quantum}
\bibinfo{author}{\bibfnamefont{M.~A.} \bibnamefont{Nielsen}} \bibnamefont{and}
  \bibinfo{author}{\bibfnamefont{I.~L.} \bibnamefont{Chuang}},
  \emph{\bibinfo{title}{Quantum computation and quantum information}}
  (\bibinfo{publisher}{Cambridge university press}, \bibinfo{year}{2010}).

\end{thebibliography}

\end{document}